\documentclass[runningheads]{llncs}

\usepackage{xargs}                      
\usepackage[pdftex,dvipsnames]{xcolor}  
\usepackage[colorinlistoftodos,prependcaption,textsize=scriptsize]{todonotes}

\newcommandx{\commentAS}[2][1=]{\todo[linecolor=blue,  backgroundcolor=blue!25,  bordercolor=blue,#1] 
            {\textbf{\sf AS:} #2}}
\newcommandx{\commentJV}[2][1=]{\todo[linecolor=Plum,  backgroundcolor=Plum!25,  bordercolor=Plum,#1] 
            {\textbf{\sf JV:} #2}}
\newcommandx{\commentGC}[2][1=]{\todo[linecolor=green, backgroundcolor=green!25, bordercolor=green,#1]
            {\textbf{\sf GC:} #2}}
\newcommandx{\commentKB}[2][1=]{\todo[linecolor=red,   backgroundcolor=red!25,   bordercolor=red,#1]  
            {\textbf{\sf KB:} #2}}

\usepackage[T1]{fontenc}
\usepackage{graphicx}
\usepackage{amsmath}
\usepackage{amssymb}
\usepackage{url}
\usepackage{booktabs}
\usepackage{array}

\usepackage{microtype}
\microtypesetup{activate=true}

\usepackage{svg}
\usepackage{tikz}
\usetikzlibrary{decorations.markings}

\DeclareMathOperator{\len}{length}

\begin{document}

\title{Maniplexes as a Foundation for Cross-Linked Databases of Symmetric Objects}

\titlerunning{Cross-Linked Databases of Symmetric Objects}

\author{%
    Katja Berčič \inst{1}\orcidID{0000-0002-6678-8975} \and
    Gabe Cunningham\inst{2}\orcidID{0000-0001-7322-6826} \and
    Andrés David Santamaría-Galvis \inst{3}\orcidID{0000-0001-5881-2886} \and
    Jano\v{s} Vidali\inst{1,3}\orcidID{0000-0001-8061-9169}%
}

\authorrunning{K. Ber\v{c}i\v{c} et al.}

\institute{%
    Faculty of Mathematics and Physics, University of Ljubljana, Slovenia \and
    School of Computing and Data Science, Wentworth Institute of Technology, Boston, MA, USA \and
    Institute of Mathematics, Physics and Mechanics, Ljubljana, Slovenia%
}

\maketitle

\begin{abstract}

    Graphs, maps on surfaces, and abstract polytopes are related combinatorial structures that tend to be studied by different communities using their own tools and databases. 
    Maniplexes provide a unifying framework that captures all of them. 
    A single database built around maniplexes would help researchers recognize shared structures and translate results across fields. 
    Here we present a compact, interoperable format for storing maniplexes as edge-labeled graphs, designed with such a database in mind. 
    As a first step, we connect two existing datasets of regular $4$-maniplexes to the House of Graphs and to Potočnik's tetravalent graph censuses, using canonical forms of their flag graphs, $1$-skeleton graphs, and $1$-coskeleton graphs.


\end{abstract}


\section{Introduction}

Graphs, maps on surfaces, abstract polytopes, and their generalizations are discrete structures studied across combinatorics, group theory, and geometry.
They are closely interconnected: every polytope of rank~$3$ is also a map, and every abstract polytope has an associated flag graph.
Yet the tools, terminologies, and databases for these objects remain largely independent, making it harder to spot connections and translate results across fields.
Maniplexes, introduced by Wilson~\cite{Wilson2012}, express maps, abstract polytopes, and related objects in a common language of edge-labeled graphs with additional structure.
Regular maniplexes, those whose color-preserving automorphism group acts regularly on vertices, are well studied (see~\cite{faithful,VoltOps,faithful-and-thin} for recent work).

Systematic study through databases requires a storage format that supports isomorphism checking, so that duplicates can be recognized and entries in different databases reliably linked.
Canonical labeling provides this: two graphs are isomorphic if and only if their canonical forms agree, reducing isomorphism testing to string comparison.
The de facto standard for graphs in databases, McKay's \texttt{graph6}/\texttt{sparse6} format~\cite{McKay2014}, is built on this principle but does not support edge labels.
The \texttt{lsparse6} format~\cite{DiscreteZOOlsparse6}, designed for the DiscreteZOO project~\cite{DiscreteZOO,bv20} which we describe here for the first time, extends \texttt{sparse6} to edge-labeled multigraphs and supports canonical forms in the same way.

\paragraph{Contributions.}
This paper makes three contributions:
\begin{enumerate}
    \item
        A self-contained description of the \texttt{lsparse6} format (Section~\ref{sec:lsparse6}) and its suitability for maniplexes.
    \item
        Classification of the contents of two datasets~\cite{maniplex-data} of regular $4$-maniplexes (Table~\ref{tab:r4m}) as polytopal, faithful, or unfaithful.    
    \item
        A cross-linked dataset (Section~\ref{sec:dataset}) associating these maniplexes with entries in the House of Graphs~\cite{HouseOfGraphs} and Poto\v{c}nik's censuses of highly symmetric tetravalent graphs~\cite{Census2AT,Graphsym,CensusAT-cubic,CensusAT-bounding,WilsonCensus,WilsonPotocnikRecipes}.%
\end{enumerate}

The remainder of the paper is organized as follows.
Section~\ref{sec:objects} introduces the relevant mathematical objects.
Section~\ref{sec:tools} describes computational tools and surveys existing databases.
Section~\ref{sec:formats} reviews the \texttt{graph6} and \texttt{sparse6} formats and introduces the \texttt{lsparse6} format for edge-labeled graphs.
Section~\ref{sec:dataset} presents the maniplex datasets and their cross-links to existing graph databases.
Section~\ref{sec:conclusions} concludes with directions for future work.

\section{Discrete and Finite Symmetric Objects}
\label{sec:objects}

We give only the essential definitions for this paper; for a detailed treatment, see~\cite{faithful,nedela2007maps,ACP}.

A \emph{graph} $G = (V, E)$ consists of a finite set of vertices $V$ and a set of edges $E \subseteq \binom{V}{2}$.
A graph is \emph{$k$-regular} (or \emph{$k$-valent}) if every vertex has exactly $k$ neighbors.
An \emph{arc} of a graph is an ordered pair of adjacent vertices.
An \emph{edge-labeled graph} (or \emph{edge-colored graph}) is a graph together with a function that assigns a label from some set $\Lambda$ to each edge.

The symmetries of a graph are captured by its \emph{automorphisms}: an automorphism of $G$ is a permutation of $V$ that preserves adjacency, and the collection of all automorphisms forms the \emph{automorphism group} $\mathrm{Aut}(G)$.
For edge-labeled graphs, automorphisms must additionally preserve edge labels.
A group \emph{acts transitively} on a set if for every pair of elements there is a group element mapping one to the other; informally, all elements ``look the same'' from the group's point of view.
A graph is \emph{vertex-transitive} if $\mathrm{Aut}(G)$ acts transitively on $V$, and \emph{arc-transitive} if it acts transitively on ordered pairs of adjacent vertices.
\begin{figure}
  \centering
  \begin{tikzpicture}[vertex/.style={circle, draw}, scale=0.7, transform shape,
    arc arrow/.style={color=blue, postaction={decorate, decoration={markings,
      mark=at position 0.25 with {\arrow{>}},
      mark=at position 0.75 with {\arrow{<}}}}}]
  \tikzset{vertex/.style={circle, draw, minimum size=6pt, inner sep=0pt}}
  \foreach \X[count=\Y] in {black,black,black,black,black}
  {\node[vertex, fill=\X] (x-\Y) at ({72*\Y+30}:2){}; }
  \draw[arc arrow, line width=0.25mm] (x-1) -- (x-2);
  \draw[arc arrow, line width=0.25mm] (x-2) -- (x-3);
  \draw[arc arrow, line width=0.25mm] (x-3) -- (x-4);
  \draw[arc arrow, line width=0.25mm] (x-4) -- (x-5);
  \draw[arc arrow, line width=0.25mm] (x-5) -- (x-1);
  \path (x-3) -- (x-4) coordinate[midway] (mid);
  \draw[dashed, shorten >=-0.3cm, shorten <=-0.3cm] (x-1) -- (mid);
  \coordinate (a1) at ({72*5+30}:1);
  \coordinate (a2) at ({72*2+30}:1);
  \draw[<->, thick] (a1) to[bend right=30] (a2);
  \draw[->, thick, shorten >=0.4cm, shorten <=0.4cm] (x-4) to[bend right=30] (x-5);
  \end{tikzpicture}
  \caption{The cycle graph $C_5$ with its arcs represented as arrows pointing away from the source vertex. 
           The dihedral group generated by the rotation by one vertex and a reflection acts transitively on vertices and arcs.}
  \label{fig:cycle}
\end{figure}
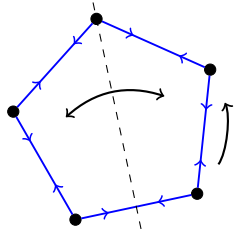
For a simple example, consider the cycle graph $C_n$ (see Figure~\ref{fig:cycle}).
Its automorphism group is the dihedral group, generated for example by the rotation by one vertex and a reflection, and it acts transitively on both vertices and arcs.
For a less obvious example, the reader may enjoy the proof without words that the automorphism group of the Petersen graph is $S_5$, acting transitively on vertices and arcs~\cite{PetersenGraph}.

Vertex-transitivity is guaranteed by construction for an important class of graphs.
A \emph{Cayley graph} $\mathrm{Cay}(\Gamma, S)$ is defined from a group $\Gamma$ and a symmetric generating set $S = S^{-1}$ not containing the identity: its vertex set is $\Gamma$, with an edge $\{g, gs\}$ for each $g \in \Gamma$ and $s \in S$.
Since $\Gamma$ acts on itself by left multiplication, every Cayley graph is vertex-transitive.
\begin{figure}
  \centering
  \includegraphics[scale=0.7]{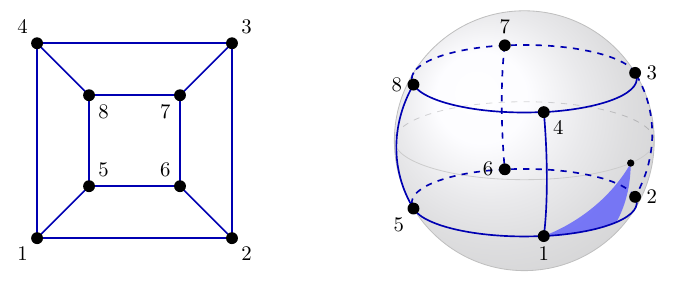}
  \caption{The cube graph and its embedding on a sphere. 
           The highlighted triangle shows the flag $(1, 12, 1234)$.}
  \label{fig:cube-embedding}
\end{figure}
A \emph{map} is a cellular embedding of a graph on a closed surface such that each face is homeomorphic to an open disk.
Equivalently, a map can be described combinatorially by its set of \emph{flags}---triples $(v, e, f)$ of a mutually incident vertex $v$, edge $e$, and face $f$---together with three fixed-point-free involutions $r_0, r_1, r_2$ acting on flags, where $r_i$ replaces the $i$-th element of the triple (see Figure~\ref{fig:cube-embedding}).
For example, $r_0$ sends a flag to the unique adjacent flag that differs in the vertex.
The composition $r_0 r_2$ has order~$2$.
A map is \emph{regular} if its automorphism group acts regularly (that is, transitively and freely) on flags.

The combinatorial description above leads naturally to a graph encoding.
The \emph{flag graph} of a map has the flags as vertices, with an edge $\{\Phi, r_i(\Phi)\}$ colored $i$ for each flag $\Phi$ and each $i \in \{0, 1, 2\}$.
Since each involution $r_i$ is fixed-point-free, every color class is a perfect matching, making the flag graph $3$-regular.
The condition that $r_0 r_2$ has order~$2$ further implies that the edges of colors $0$ and $2$ form $4$-cycles.
A map is \emph{orientable} if the underlying surface is orientable (genus $g \geq 0$); otherwise it is \emph{non-orientable}.


An \emph{abstract polytope} of rank~$n$ is a partially ordered set $(\mathcal{P}, \leq)$ with a rank function $\mathrm{rk}\colon \mathcal{P} \to \{-1, 0, 1, \dotsc, n\}$ satisfying several axioms: there is a unique element of rank~$-1$ and a unique element of rank~$n$, every maximal chain contains $n+2$ elements, and the poset is strongly connected and satisfies the diamond property~\cite{McMullenSchulte2002}.
For our purposes, it will suffice to think of the diamond property as generalizing the pattern of two edges meeting at a vertex and two faces meeting at an edge in a convex polyhedron.
Elements of rank~$i$ are called \emph{$i$-faces} of the polytope. 

\begin{figure}
  \centering
  \includegraphics[scale=0.44]{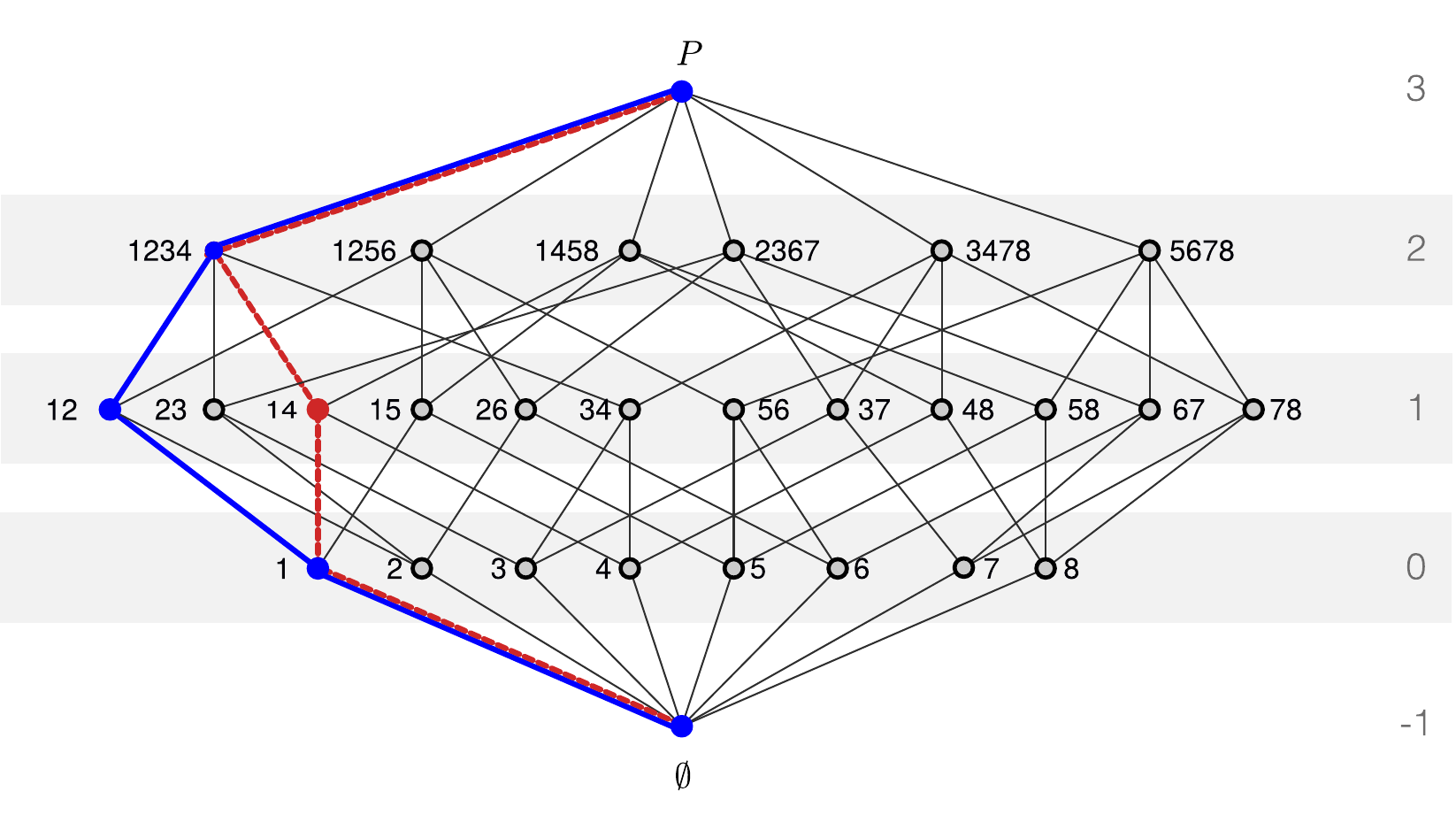}
  \caption{The cube as an abstract polytope.
           The blue (solid) flag, $(\emptyset, 1, 12, 1234, P)$, is adjacent to the red (dashed) flag $(\emptyset, 1, 14, 1234, P)$.
           The ranks of facets are shown on the right for convenience.}
  \label{fig:poset}
\end{figure}

As with maps, the combinatorics of an abstract polytope can be encoded via flags and adjacency.
A \emph{flag} of $\mathcal{P}$ is a maximal chain; two flags are \emph{$i$-adjacent} if they differ in exactly the element of rank~$i$ (see Figure~\ref{fig:poset}).%
\begin{figure}
  \centering
  \includegraphics[scale=0.4]{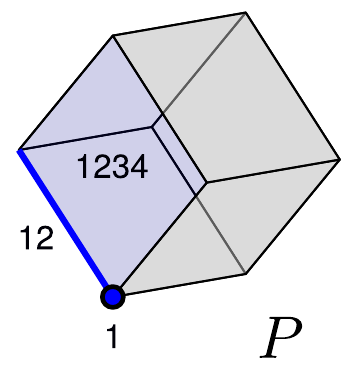}
  \caption{The cube with faces of the flag $(\emptyset, 1, 12, 1234, P)$ highlighted.}
  \label{fig:flag}
\end{figure}
Each flag has exactly one $i$-adjacent flag for each $i \in \{0, \dotsc, n{-}1\}$, giving rise to $n$ fixed-point-free involutions $r_0, \dotsc, r_{n-1}$ on flags.
These involutions generate the \emph{monodromy group}, and the resulting flag graph is defined exactly as for maps, now with $n$ colors.
An abstract polytope is \emph{regular} if its automorphism group acts regularly on flags.
The automorphism group of a regular polytope is generated by elements $\rho_0, \dotsc, \rho_{n-1}$, where each $\rho_i$ sends a fixed base flag $\Phi$ to its $i$-adjacent flag $\Phi^i$; these generators satisfy the intersection condition
\begin{equation}
    \langle \rho_i : i \in I \rangle 
      \cap 
    \langle \rho_j : j \in J \rangle =
    \langle \rho_k : k \in I \cap J \rangle 
    \ 
    \text{ for all } 
    \ 
    I, J \subseteq \{0,\dotsc,n{-}1\}.
    \label{eq:int-cond}
\end{equation}

Maniplexes~\cite{Wilson2012} generalize both maps and abstract polytopes, capturing both within a single graph-theoretic framework. We adopt the flag graph definition, which has become standard.

\begin{definition}
    A \emph{maniplex} of rank~$n$ is a connected $n$-regular graph $M$ with a
    proper edge-coloring by colors $\{0, 1, \dotsc, n{-}1\}$ such that:
    \begin{enumerate}
        \item
            For each color $i$, the edges of color $i$ form a perfect matching.
        \item
            For any two colors $i,j$ with $|i-j| \geq 2$, the composition $r_i r_j$ has order~$2$, where $r_i$ denotes the involution sending each vertex to its unique $i$-colored neighbor.
    \end{enumerate}
\end{definition}
\begin{figure}
    \centering
    \includegraphics[width=0.53\linewidth]{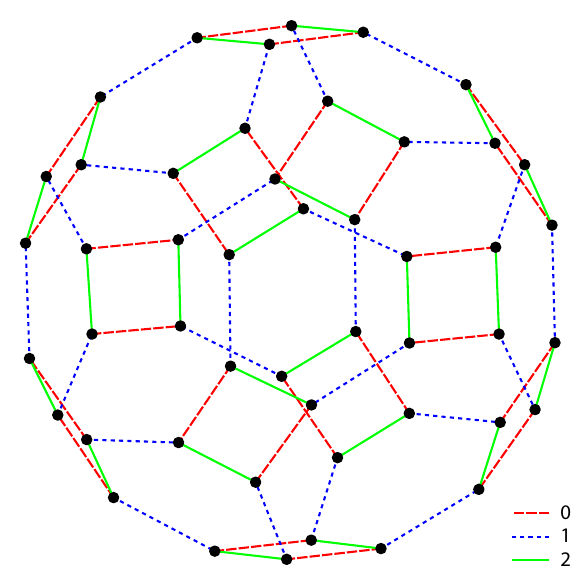}
    \caption{The flag graph of the cube. 
             The edges of colors $0$, $1$, and $2$ are shown in red (dashed), blue (dotted), and green (solid), respectively. The edges of colors $0$ and $2$ form $4$-cycles, corresponding to edges in the cube; the edges of colors $1$ and $2$ form $6$-cycles, corresponding to the vertices of the cube; and the edges of colors $0$ and $1$ form $8$-cycles, corresponding to the faces of the cube.}
    \label{fig:maniplex} 
\end{figure}  

The vertices of a maniplex are called \emph{flags}, and a maniplex of rank~$n$ is its own flag graph.
A maniplex is \emph{regular} if its group of color-preserving automorphisms acts regularly on vertices.
For a regular maniplex, this automorphism group is generated by involutions $\rho_0, \dotsc, \rho_{n-1}$; in fact, there is an isomorphism from the automorphism group to the monodromy group sending each $\rho_i$ to $r_i$, although the two groups act on flags differently.

We will also need some derived structures:
The \emph{underlying graph of the flag graph} is the flag graph without its edge coloring.
The \emph{$i$-faces} of a maniplex $M$ of rank~$n$ are the connected components of the subgraph of $M$ obtained by deleting all edges of color $i$; for a polytopal maniplex, these agree with the $i$-faces of the underlying polytope.
The \emph{$1$-skeleton} of $M$ is the graph whose vertices are the $0$-faces of $M$ and whose edges are the $1$-faces, with incidence induced from $M$.
The \emph{$1$-coskeleton} of $M$ is the graph formed by the $1$-skeleton of the dual, i.e., its vertices are the $(n - 1)$-faces, and the edges are the $(n - 2)$-faces of $M$.

Every abstract polytope of rank~$n$ yields a maniplex of rank~$n$ via its flag graph, but not every maniplex is \emph{polytopal} (isomorphic to the flag graph of an abstract polytope).
The class of regular $4$-maniplexes can be partitioned into three subclasses:
\begin{itemize}
  \item 
    \textbf{Polytopal (P):} the maniplex is isomorphic to the flag graph of an abstract polytope. 
    (In the case of regular maniplexes, this means that the intersection condition from Equation~\ref{eq:int-cond} holds).
  \item 
    \textbf{Faithful (F):} the maniplex is faithful but not polytopal.
    Every polytopal maniplex is faithful, meaning it can be represented by a partially ordered set analogous to an abstract polytope without loss of information~\cite{faithful}.
  \item 
    \textbf{Unfaithful (U):} the maniplex is unfaithful, so we cannot recover the maniplex from the partially ordered set of its faces under incidence.
\end{itemize}

Orientability for maniplexes generalizes that for maps: a maniplex $M$ of rank~$n$ is \emph{orientable} if it is a bipartite graph~\cite{Wilson2012}.

\section{Computational Tools and Existing Data}
\label{sec:tools}

\subsection{Automorphism Groups and the RAMP Package}
\label{sec:ramp}

The automorphism group of a regular maniplex of rank~$n$ is a group $\Gamma$ generated by $n$ involutions $\rho_0, \dotsc, \rho_{n-1}$ satisfying the relations $(\rho_i \rho_j)^2 = 1$ whenever $|i - j| \geq 2$.
Conversely, given such a group with generating set $\{\rho_0, \dotsc, \rho_{n-1} \}$, one recovers the regular maniplex as the Cayley graph of $\Gamma$ with these generators.
For rank~$n=4$, this Cayley graph is tetravalent, since each generator is an involution.

The \textbf{RAMP} (Research Assistant for Maniplexes and Polytopes) package~\cite{RAMP} for the computer algebra system GAP~\cite{GAP} provides tools for constructing and analyzing maniplexes, maps, and abstract polytopes from their automorphism or monodromy groups.
RAMP can compute the flag graph of a maniplex from its monodromy group, yielding the edge-labeled graph that we then encode in \texttt{lsparse6} format.

\subsection{A Brief Survey of Symmetric Object Data}
\label{sec:existing}

Several databases of symmetric discrete structures are relevant to this work.
We survey them here, organized by the type of object they store.

\paragraph{Graphs.}
The House of Graphs~\cite{HouseOfGraphs} is a database of ``interesting'' graphs.
It supports efficient lookup by canonical form, allows users to add new graphs, and stores computed invariants alongside each entry.
Poto\v{c}nik and collaborators have compiled censuses of highly symmetric tetravalent graphs available at Graphsym~\cite{Graphsym}, including arc-transitive graphs on up to $640$ vertices~\cite{CensusAT-cubic,CensusAT-bounding} (with Spiga and Verret), $2$-arc-transitive graphs on up to $2000$ vertices~\cite{Census2AT}, and edge-transitive graphs on up to $512$ vertices~\cite{WilsonCensus,WilsonPotocnikRecipes} (with Wilson).
The Encyclopedia of Graphs~\cite{EncyclopediaOfGraphs} is an atlas of named and notable graphs with precomputed properties.
DistanceRegular.org~\cite{DistanceRegular} is a repository of distance-regular graphs, a class that overlaps with highly symmetric graphs arising from polytopes.

\paragraph{Maps and hypermaps.}
Conder maintains extensive lists of regular and chiral maps, hypermaps, and polytopes for small groups~\cite{ConderLists}, a key reference for maps on surfaces.

\paragraph{Abstract polytopes.}
Leemans and collaborators provide an atlas of abstract regular polytopes for small almost simple groups~\cite{LeemansAtlas}.
Hartley's Atlas of Small Regular Polytopes~\cite{HartleyAtlas} and the companion Atlas of Small Chiral Polytopes~\cite{HartleyChiral} enumerate polytopes by group order.
A related dataset covers polytopes derived from sporadic simple groups~\cite{HartleySporadic}.

\paragraph{Groups.}
The Small Groups Library~\cite{SmallGroups} enumerates all finite groups up to isomorphism for small orders, and is distributed with GAP~\cite{GAP} and Magma.

\paragraph{Cross-cutting projects.}
DiscreteZOO~\cite{DiscreteZOO} is an ongoing project to compile and cross-link datasets of discrete combinatorial objects.
It introduced the \texttt{lsparse6} format for compact storage of edge-labeled graphs.

No existing database targets maniplexes specifically, and no standard format for storing them existed before DiscreteZOO introduced \texttt{lsparse6}.

\section{Graph Formats}
\label{sec:formats}

McKay's \texttt{graph6} and \texttt{sparse6} formats~\cite{GraphFormats} are standard formats to represent an undirected graph as an ASCII string.
In both formats, the graph has vertex set $\{0, 1, \dotsc, n - 1\}$ and it is represented as a sequence of printable ASCII bytes (values 63--126), enabling storage, transmission, and processing with standard tools.
The format \textbf{\texttt{graph6}} encodes the upper triangle of the adjacency matrix as a bit string.  
For a graph on $n$ vertices, it uses $\lceil n(n-1)/12 \rceil + O(\log n)$ bytes, making it suitable for dense graphs.
The format \textbf{\texttt{sparse6}} encodes only the edge list.
It encodes edges into a bitstream using pairs of a flag bit and an incoming endpoint (encoded in $\lceil \log_2 n \rceil$ bits).
The numerically larger of the two endpoints is tracked implicitly and changes based on the flag bit by up to~$1$, or is given explicitly if the incoming endpoint is larger than it.
If the incoming endpoint is smaller, then the tracked and incoming endpoints are connected.
Bits are packed into $6$-bit groups (big-endian, padded with $1$-bits) and stored as bytes in the range $[63, 126]$.
For $k$-regular graphs on $n$ vertices, the encoding uses $O(kn \log n)$ bytes, which is significantly smaller than \texttt{graph6} for sparse or moderate degree graphs.%

Both formats begin with a short header that identifies the format (e.g., \texttt{:} for \texttt{sparse6}) and encodes $n$ using McKay's variable-length $N$ function, which encodes smaller integer values with fewer bytes.
\[
  \len (N(n)) = \begin{cases}
    \text{one byte}         & 0 \leq n \leq 62, \\
    \text{four bytes}       & 63 \leq n \leq 258047, \\
    \text{eight bytes}      & 258048 \leq n \leq 68719476735.
  \end{cases}
\]

A \emph{canonical form} of a graph $G$ is a graph $\mathrm{can}(G) \cong G$ chosen by a deterministic rule so that $\mathrm{can}(G) = \mathrm{can}(H)$ if and only if $G \cong H$.  The \texttt{nauty/Traces} software~\cite{McKay2014} computes canonical labelings of simple graphs efficiently, even for large graphs.
The canonical \texttt{graph6} or \texttt{sparse6} string of a graph then serves as a unique identifier suitable for database lookup and deduplication. 
Efficient canonical-form computation makes \texttt{(l)sparse6}-based formats well suited for cross-referencing:
two independently produced encodings of isomorphic graphs will yield the same canonical string, regardless of vertex ordering.
Canonical labeling only relabels vertices; edge labels (and thus maniplex colors) are preserved.

\subsection{The \texttt{lsparse6} Format}
\label{sec:lsparse6}

The \texttt{lsparse6} format~\cite{DiscreteZOOlsparse6} extends \texttt{sparse6} to undirected edge-labeled multigraphs.
An \texttt{lsparse6} string for a graph with $n$ vertices, $m$ edges, and $\ell$ distinct edge labels consists of three concatenated parts: the \texttt{sparse6} encoding of the underlying unlabeled graph as specified by McKay~\cite{McKay2014}, a \texttt{\#} delimiter (ASCII code 35), and the label data.

The label data is encoded as follows.
The number of distinct labels $\ell$ is written first using McKay's $N$ function.
Writing $k = \lceil \log_2 \ell \rceil$, each edge label (an integer in $\{0, \dotsc, \ell{-}1\}$) is encoded in $k$ bits, listed in the same order as the corresponding edges in the \texttt{sparse6} encoding.
The bits are packed into $6$-bit groups in big-endian order, padded with $1$-bits to reach a multiple of $6$, and stored as bytes in the range $[63, 126]$ by adding $63$.

For a regular $4$-maniplex with $f$ flags, the flag graph has $f$ vertices and $2f$ edges (each vertex is incident to $4$ edges, giving $4f/2 = 2f$ undirected edges).
The edge labels lie in $\{0,1,2,3\}$, so $\ell = 4$ and $k = 2$ bits per label.
The total size of the \texttt{lsparse6} string is therefore the size of the underlying \texttt{sparse6} string plus one byte for the label count plus $\lceil 2 \cdot 2f / 6 \rceil$ additional bytes for the labels---a modest overhead of roughly $2f/3$ bytes beyond the graph structure, which itself takes $O(f \log f)$ bytes in \texttt{sparse6}.



\paragraph{Relation to \tt{\tt sparse6}.}
When $\ell = 1$ (all edges have the same label), the label data reduces to a
single byte (value $64 = \texttt{`@'}$) 
and the resulting string is only marginally longer than the underlying \texttt{sparse6} string.
(Readers with unlabeled graphs should, of course, simply use \texttt{sparse6}.)

\subsection{Multi-edges and Loops}

The maniplex community frequently works with multigraphs, that is, graphs that may have multiple edges between the same pair of vertices as well as loops (edges from a vertex to itself).
While the present paper does not require them, we briefly discuss support for these features with an eye toward future extensions of our tools.
Since \texttt{lsparse6} is built on \texttt{sparse6}, it inherits support for both multi-edges and loops.
Table~\ref{tab:graph-formats} compares the compact graph representations discussed so far.
\begin{table}[h]
  \centering
  \caption{Comparison of compact graph formats. Support for multi-edges and loops is inherited from \texttt{sparse6}.}
  \label{tab:graph-formats}
  \begin{tabular}{p{0.15\textwidth} p{0.35\textwidth} p{0.15\textwidth} p{0.15\textwidth} p{0.15\textwidth}}
    \toprule
    Format & Best for & Multi-edges & Loops & Edge labels \\
    \midrule
    \texttt{graph6}   & Small or dense graphs      & No        & No        & No  \\
    \texttt{sparse6}  & Large sparse graphs        & Yes       & Yes       & No  \\
    \texttt{lsparse6} & Edge-labeled (multi)graphs & Inherited & Inherited & Yes \\
    \bottomrule
  \end{tabular}
\end{table}
For isomorphism computation and canonical labeling, both \texttt{nauty} and \texttt{Traces} are primarily designed for simple graphs~\cite{McKay2014}; Table~\ref{tab:nauty-traces} compares their capabilities.
In particular, edge-colored multigraphs are not natively supported and require workarounds
or auxiliary tools such as \texttt{multig}~\cite{McKay2014}.
\begin{table}[h]
  \centering
  \caption{Capabilities of \texttt{nauty} and \texttt{Traces}.}
  \label{tab:nauty-traces}
  \begin{tabular}{p{0.25\textwidth} p{0.28\textwidth} p{0.45\textwidth}}
    \toprule
    Feature & \texttt{nauty} & \texttt{Traces} \\
    \midrule
    Simple graphs      & $\checkmark$ (native)        & $\checkmark$ (native)              \\
    Vertex coloring    & $\checkmark$                 & $\checkmark$                       \\
    Edge coloring      & $\times$ (requires encoding) & $\times$ (requires encoding)       \\
    Multiple edges     & $\times$ (requires encoding) & $\times$ (requires encoding)       \\
    Loops              & limited                      & limited                            \\
    Canonical labeling & $\checkmark$ & $\checkmark$ (often faster on large sparse graphs) \\
    \bottomrule
  \end{tabular}
\end{table}
%

\section{A Cross-Linked Maniplex Dataset}
\label{sec:dataset}

The cross-linking pipeline has so far been completed for all maniplexes with up to $1000$ flags; the remaining entries will be added as computations finish.
An up-to-date version of the cross-linked dataset is available at~\cite{crosslinked-data}.

Our source maniplex data~\cite{maniplex-data} consist of two complete enumerations of regular $4$-maniplexes: all non-orientable ones with up to $3550$ flags and all orientable ones with up to $7100$ flags (Table~\ref{tab:r4m}), each classified as polytopal, faithful, or unfaithful.
\begin{table}[h]
  \centering
  \caption{Overview of the two datasets of regular $4$-maniplexes: all non-orientable with up to $3550$ flags and all orientable with up to $7100$ flags.}
  \label{tab:r4m}
  \begin{tabular}{l *{2}{>{\raggedleft\arraybackslash}p{0.1\textwidth} >{\raggedleft\arraybackslash}p{0.08\textwidth} >{\raggedleft\arraybackslash}p{0.1\textwidth} >{\raggedleft\arraybackslash}p{0.1\textwidth}}}
    \toprule
    & \multicolumn{4}{c}{Non-orientable} & \multicolumn{4}{c}{Orientable} \\
    \cmidrule(l){2-5} \cmidrule(l){6-9}
    Type & Count & Min. & Max. & $\leq 1000$ & Count & Min. & Max. & $\leq 1000$ \\
    \midrule
    Polytopal (P)  & 5763     & 48 & 3528 & 750      & 77\,887  & 16 & 7100 & 2778 \\
    Faithful (F)   & 10\,547  & 72 & 3528 & 958      & 345\,212 & 32 & 7088 & 5143 \\
    Unfaithful (U) & 149\,747 & 8  & 3548 & 17\,712  & 164\,695 & 8  & 7096 & 5293 \\
    \midrule
    Total          & 166\,057 &    &      & 19\,420  & 587\,794 &    &      & 13\,214 \\
    \bottomrule
  \end{tabular}
\end{table}
Each regular $4$-maniplex gives rise to three graphs that can be matched against existing datasets: the underlying (unlabeled) flag graph, which is tetravalent, and the $1$-skeleton and $1$-coskeleton (Section~\ref{sec:objects}), which are typically much smaller.
We cross-link all three against the House of Graphs~\cite{HouseOfGraphs} and three of Poto\v{c}nik's tetravalent graph censuses: the census of arc-transitive tetravalent graphs on up to $640$ vertices~\cite{CensusAT-bounding,CensusAT-cubic} (joint work with Spiga and Verret), the census of $2$-arc-transitive tetravalent graphs on up to $2000$ vertices~\cite{Census2AT}, and the census of edge-transitive tetravalent graphs on up to $512$ vertices~\cite{WilsonCensus,WilsonPotocnikRecipes} (joint work with Wilson).
The contents of these databases are summarized in Table~\ref{tab:graphs}.
\begin{table}[h]
  \centering
  \caption{Graph databases used for cross-linking.}
  \label{tab:graphs}
  \begin{tabular}{p{0.25\textwidth} >{\raggedleft\arraybackslash}p{0.15\textwidth} >{\raggedleft\arraybackslash}p{0.13\textwidth} @{\hspace{2em}} p{0.40\textwidth}}
    \toprule
    Dataset & Max order & Graphs & Comment \\
    \midrule
    House of Graphs~\cite{HouseOfGraphs}                      & 250  & 28497 & Not exhaustive \\
    Edge-transitive~\cite{WilsonCensus,WilsonPotocnikRecipes} & 512  & 7364  & Complete for several segments \\
    Arc-transitive~\cite{CensusAT-bounding,CensusAT-cubic}    & 640  & 4820  & Complete \\
    $2$-arc-transitive~\cite{Census2AT}                       & 2000 & 165   & Complete \\
    \bottomrule
  \end{tabular}
\end{table}
The three censuses have increasing symmetry requirements: $2$-arc-transitivity implies arc-transitivity, which in turn implies edge-transitivity.
Restricting to the range of orders covered by a less restrictive census, the graphs of a stricter census are contained in the less restrictive one.

\begin{figure}[ht!]
  \centering
  \includegraphics[width=\textwidth]{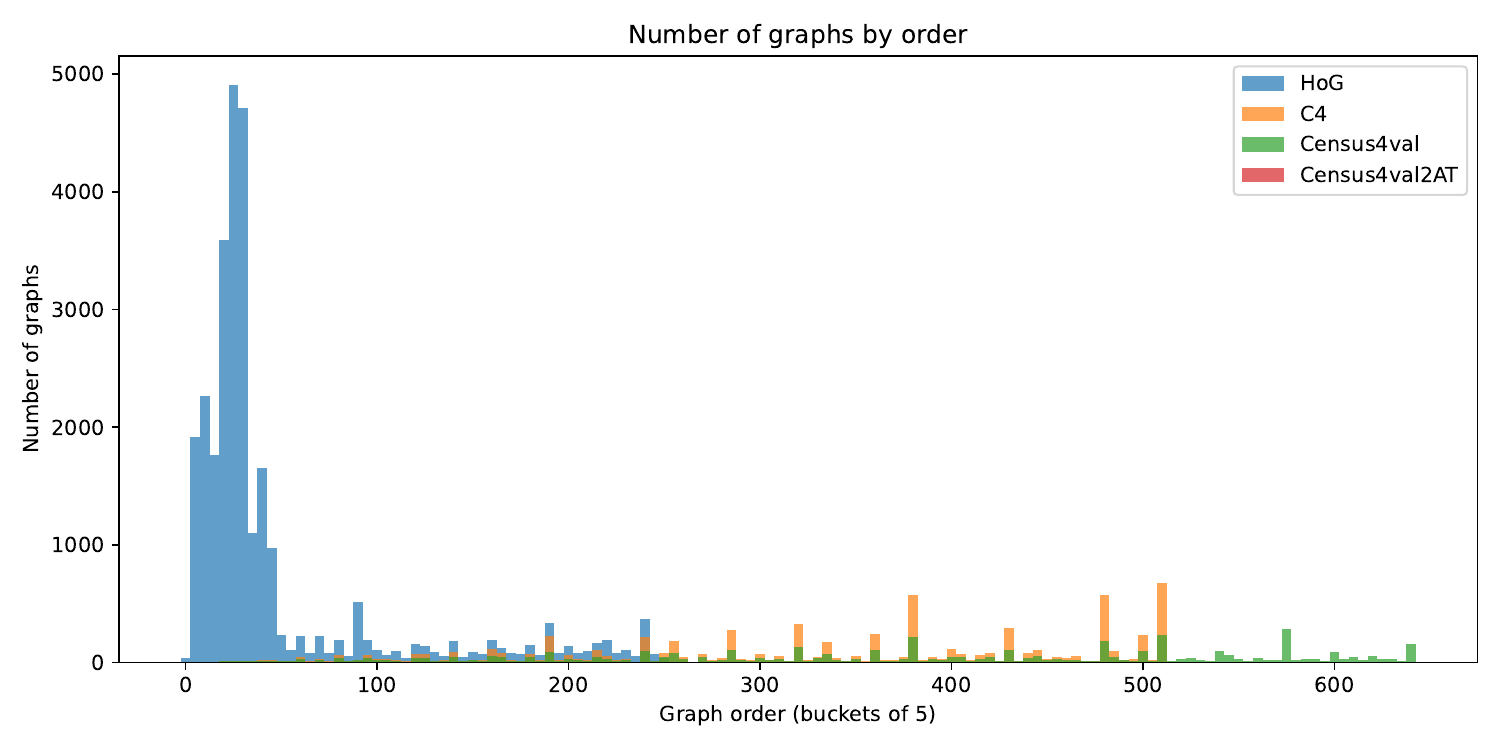}
  \caption{Distribution of graphs by order in each database, showing orders up to $640$. 
           The only dataset extending beyond this range is the $2$-arc-transitive census, where graphs are very rare. 
           Graphs are grouped into buckets of $5$ consecutive orders for better visualization.}
  \label{fig:graphs-by-order}
\end{figure}

The cross-linking proceeds via the following pipeline.
\begin{enumerate}
  \item 
    Data from~\cite{maniplex-data} are imported into RAMP and classified as polytopal, faithful, or unfaithful. (See Table~\ref{tab:r4m}.)
  \item 
    We use RAMP to compute the flag graph of each maniplex and \texttt{nauty} to compute its canonical form.
  \item 
    Each flag graph is encoded as an \texttt{lsparse6} string.
    We strip the edge labels to obtain the underlying (unlabeled) tetravalent graph.
  \item 
    We use RAMP and \texttt{nauty} to compute canonical forms of the $1$-skeleton and $1$-coskeleton.
  \item
    We compute canonical forms of the entries in the House of Graphs and the three censuses, and link each maniplex to any graph entry with a matching canonical form for the underlying flag graph, $1$-skeleton, or $1$-coskeleton.
\end{enumerate}
The pipeline currently uses the following software:
GAP~\cite{GAP} version~4.15.1,
RAMP~\cite{RAMP} version~0.7.01, and
PassageMath version~10.6.44.

\subsection{Database Description}
\label{sec:db-description}

\begin{figure}[htbp]
  \centering
  \includegraphics[scale=0.150]{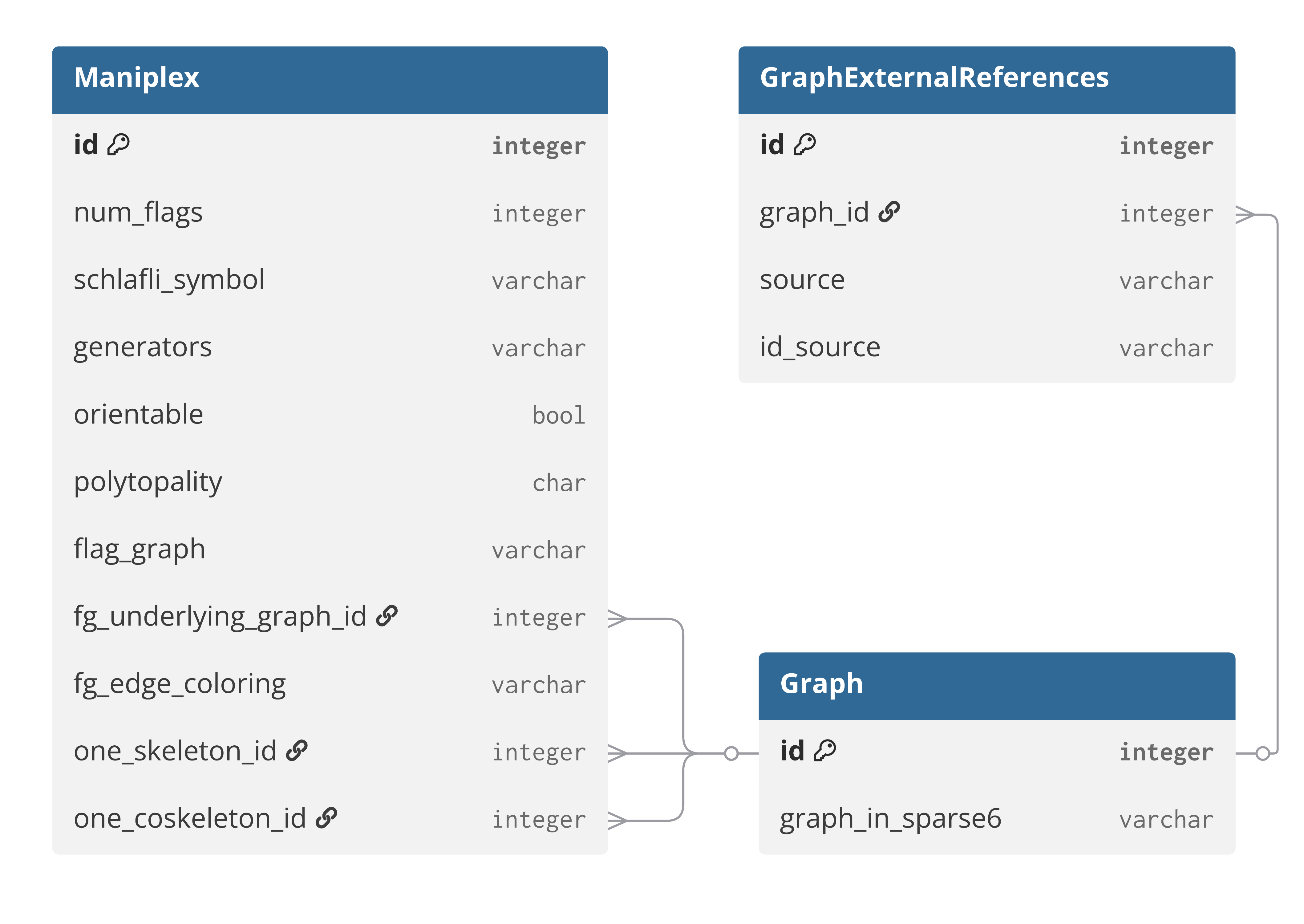}
  \caption{A schematic diagram of the cross-linked dataset.  
           Each maniplex entry (with its \texttt{lsparse6} encoding) is linked to one or more graph entries (with their \texttt{sparse6} encodings) via the canonical forms of the underlying graph, $1$-skeleton, and $1$-coskeleton. 
           The graph table also includes graph-theoretic invariants.}
    \label{fig:db_diagram}
\end{figure}

The database we use to store the cross-linking information consists of three main tables: \texttt{Maniplex}, \texttt{Graph}, and \texttt{GraphExternalReferences} (see Figure~\ref{fig:db_diagram}).
\begin{enumerate}
  \item
    The \texttt{Maniplex} table stores regular maniplexes, each identified by its number of flags, Schläfli symbol, generators, orientability, and polytopality. 
    It also stores its flag graph in \texttt{lsparse6} format.  
    Each maniplex is associated with three graph objects: the underlying graph of its flag graph, its $1$‑skeleton, and its $1$‑coskeleton.
  \item
    The \texttt{Graph} table stores simple graphs encoded in \texttt{sparse6} format. 
    These graphs appear referenced by the corresponding maniplex in the \texttt{Maniplex} table.
  \item
    The \texttt{GraphExternalReferences} table links graphs to external databases. 
    For each graph, it records the source (e.g., House of Graphs, Census4val, C4) and the identifier used in that external source.
\end{enumerate}
The table \texttt{Graph} connects to \texttt{GraphExternalReferences} (one‑to‑many), and the table \texttt{Maniplex} connects to \texttt{Graph} via three foreign keys, one for each of the three derived graph structures.

\subsection{Cross-Linking Results}
\label{sec:cross-linking-results}

We now summarize the results of the cross-linking pipeline for the $32\,634$ maniplexes with up to $1000$ flags.
Not every derived graph falls within the vertex ranges of the external datasets; Table~\ref{tab:graph-sizes} shows the coverage for flag graphs.
All $1$-skeletons and $1$-coskeletons have fewer than $250$ vertices, so they are within range of every dataset.
\begin{table}[ht]
  \centering
  \caption{Number of maniplexes (with up to $1000$ flags, $32\,634$ in total) whose flag graphs fall within the order ranges of the databases in Table~\ref{tab:graphs}.}
  \label{tab:graph-sizes}
  \begin{tabular}{l >{\raggedleft\arraybackslash}p{0.2\textwidth} >{\raggedleft\arraybackslash}p{0.2\textwidth} >{\raggedleft\arraybackslash}p{0.2\textwidth} >{\raggedleft\arraybackslash}p{0.2\textwidth}}
    \toprule
    Derived graph & $\leq 250$ & $251$--$512$ & $513$--$640$ & $> 640$ \\
    \midrule
    Flag graph     & 3192 & 8722 & 4199 & 16\,521 \\
    \bottomrule
  \end{tabular}
\end{table}
Table~\ref{tab:matches} shows the number of maniplexes (with up to $1000$ flags) for which a match was found in each database, broken down by the type of derived graph.
Different maniplexes can share the same underlying graph: a single tetravalent graph may admit several distinct edge-colorings satisfying the maniplex axioms.
The $1$-coskeleton columns are omitted from this table: the set of graphs occurring as $1$-skeletons is the same as the set occurring as $1$-coskeletons, although for a given maniplex the two generally differ (they coincide in 5106 cases).
\begin{table}[h]
  \centering
  \caption{Number of maniplexes (with up to $1000$ flags) with at least one match in each database, by derived graph type, and the number of distinct graphs matched.}
  \label{tab:matches}
  \begin{tabular}{l *{2}{>{\raggedleft\arraybackslash}p{0.185\textwidth}} *{2}{>{\raggedleft\arraybackslash}p{0.2\textwidth}}}
    \toprule
    & \multicolumn{2}{c}{Flag graph} & \multicolumn{2}{c}{$1$-skeleton} \\
    \cmidrule(l){2-3} \cmidrule(l){4-5}
    Database & maniplexes & distinct graphs & maniplexes & distinct graphs \\
    \midrule
    House of Graphs          & 22 & 10 & 32126 & 202 \\
    Edge-transitive          & 20 & 18 & 2144  & 51  \\
    Arc-transitive           & 21 & 19 & 2144  & 51  \\
    $2$-arc-transitive       &  4 &  2 & 852   & 9   \\
    \midrule
    Total (distinct)         & 34 & 22 & 32126 & 202 \\
    \bottomrule
  \end{tabular}
\end{table}
For $1$-skeletons, all $9$ $2$-arc-transitive graphs are already among the $51$ arc-transitive ones, which are the same $51$ graphs as in the edge-transitive census; all $51$ are also included in the House of Graphs.
For flag graphs, the overlap is similar but not identical: of the $10$ flag graphs found in the House of Graphs, $7$ also appear in at least one census, while $3$ are only in the House of Graphs. One flag graph appears in the arc-transitive census but not in the edge-transitive census, as it exceeds the latter's $512$-vertex limit.
In total, there are $22$ distinct underlying graphs of flag graphs and $218$ distinct $1$-skeleton graphs with external references.
There are few matches for flag graphs, while on the other hand, only $508$ maniplexes have no skeleton match in any of the databases.

For flag graphs, the symmetry conditions on maniplexes and on the graphs in the censuses are of a different nature:
a regular maniplex has a group acting regularly on its flags, which makes the flag graph vertex-transitive, but does not in general force it to be arc- or edge-transitive in the graph-theoretic sense; conversely, a tetravalent graph in one of the censuses need not admit a proper edge-coloring satisfying the maniplex axioms.
We expect to see significantly more links as additional types of objects, such as maps, hypermaps, and abstract polytopes, are incorporated into the dataset.

\section{Conclusions and Future Work}
\label{sec:conclusions}

We have described the \texttt{lsparse6} format for encoding edge-labeled graphs and argued for its suitability as a storage format for maniplexes.
Using the RAMP package for GAP, we have encoded the maniplexes from the source datasets in (canonically labeled) \texttt{lsparse6}.
We have described a cross-linked dataset connecting these maniplexes to the House of Graphs and to censuses of highly symmetric tetravalent graphs, using canonical forms for matching.

Directions for future work include cross-linking to further databases, such as the Small Groups Library~\cite{SmallGroups} for the automorphism groups of the maniplexes, as well as map and abstract polytope datasets.


\begin{credits}
    \subsubsection{\ackname}
    
    Ber\v{c}i\v{c} was supported by the Air Force Office of Scientific Research (award number FA9550-21-1-0024), and by the Slovenian Research and Innovation Agency (program no.\ P1-0294).
    Vidali was supported by the Slovenian Research and Innovation Agency (program no.\ P1-0285).
    Ber\v{c}i\v{c}, Santamar\'{\i}a-Galvis, and Vidali were supported by the AI For Math Fund, a program of Renaissance Philanthropy.
    Cunningham was supported by the Slovenian Research and Innovation Agency (BI-US/24-26-025).
\end{credits}

\bibliographystyle{splncs04}
\bibliography{references}

\end{document}